\documentclass[titlepage,12pt]{article}

\usepackage{graphicx}
\usepackage{amssymb}

\newcommand{\R}{\mathbb R}
\newcommand{\C}{\mathbb C}

\begin{document}

\baselineskip=18pt

\begin{center}
{\Large{\bf Stability and Hopf Bifurcation in the
 \\ Watt Governor System}}
\end{center}

\begin{center}
{\bf Jorge Sotomayor, Luis Fernando Mello and \\ Denis de Carvalho
Braga}
\end{center}

\begin{center}
{\bf Abstract}
\end{center}

\vspace{0.1cm}

In this paper we study the Lyapunov stability and Hopf bifurcation
in a system coupling a Watt-centrifugal-governor with a
steam-engine. Sufficient conditions for the stability of the
equilibrium state in terms of the physical parameters and of the
bifurcating periodic orbit at most critical parameters on the
bifurcation surface are given.

\vspace{0.1cm}



\section{\bf Introduction}\label{intro}

The Watt centrifugal governor is a device that automatically
controls the speed of an engine. Dating to 1788, it can be taken
as the starting point for  automatic control theory (see
MacFarlane \cite{mac} and references therein). In this paper the
system coupling the Watt-centrifugal-governor and the steam-engine
will be called simply the Watt Governor System (WGS).

Landmarks for the study of the local stability of the WGS are the
works of Maxwell \cite{max} and Vyshnegradskii \cite{vysh}. A
simplified version of the WGS local stability based on the work of
Vyshnegradskii  is presented by Pontryagin \cite{pon}. A local
stability study  generalized to a more general Watt governor
design was carried out by  Denny \cite{denny}.

Enlightening historical comments about the Watt governor local
mathematical stability and oscillatory analysis can be found in
MacFarlane \cite{mac} and Denny \cite{denny}. There, as well as in
\cite{pon}, we learn that toward the mid $XIX$ century,
improvements in the engineering design led to  less reliable
operations in the WGS, leading to fluctuations and oscillations
instead of the  ideal stable constant speed  output requirement.
The first mathematical analysis of the stability conditions and
subsequent indication of the modification in design to avoid  the
problem was due to Maxwell \cite{max} and, in a user friendly
style likely to be better understood by engineers, by
Vyshnegradskii \cite{vysh}.

From the mathematical point of view, the  oscillatory, small
amplitude, behavior in the WGS can be  associated to a periodic
orbit that appears from a Hopf bifurcation. This was established
by Hassard et al. in \cite{has1} and  Al-Humadi and Kazarinoff in
\cite{humadi}. Another procedure, based in the method of harmonic
balance,  has been suggested by Denny \cite{denny}  to detect
large amplitude oscillations.

In this paper we present a simple way to understand the Hopf
bifurcation in a WGS, which is more general  than that presented
by Pontryagin \cite{pon}, Al-Humadi and Kazarinoff  \cite{humadi}
and Denny \cite{denny}. We believe that our approach has the
advantage of connecting the physical parameters of the system to
the stability of the stationary equilibrium point and the
bifurcating periodic orbit. It also permits a neater geometric
synthesis of the bifurcation analysis based on the algebraic
expression and geometric location of the curve ---the codimension
2 Hopf points--- characterizing the stability versus the
instability at the critical weak focal equilibria on the critical
surface of parameter values ---the codimension 1 Hopf points---
leading to the bifurcation of periodic orbits, whose stability
depend on the side of the curve at which the parameters cross the
surface.

This paper is organized as follows. In Section \ref{watt} we
introduce the differential equations that model a general WGS .
The stability of the equilibrium point of this model is analyzed
and a general version of the stability condition is obtained  and
presented in the terminology of  Vyshnegradskii.  The codimension
1 Hopf bifurcation for the WGS differential equations is studied
in Section \ref{hopf}. An  expression  --- neater than that found
in the current literature --- which determines the sign of the
first Lyapunov coefficient is obtained. Sufficient conditions for
the stability for the bifurcating periodic orbit are given.
Concluding comments are presented in Section \ref{conclusion}.

\section{A  general Watt governor system}\label{watt}

\newtheorem{teo}{Theorem}[section]
\newtheorem{lema}[teo]{Lemma}
\newtheorem{prop}[teo]{Proposition}
\newtheorem{cor}[teo]{Corollary}
\newtheorem{remark}[teo]{Remark}
\newtheorem{example}[teo]{Example}

\subsection{General differential equations}\label{diffequat}

\begin{figure}[!h]
\centerline{
\includegraphics[width=12cm]{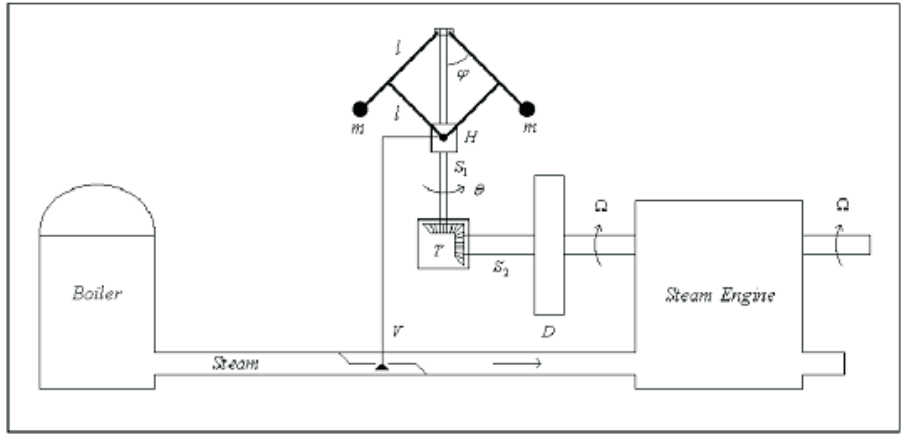}}

\caption{{\small General-Watt-centrifugal-governor-steam-engine
system}.}

\label{wattgov}

\end{figure}

The WGS studied in this paper is shown in Fig. \ref{wattgov}.
There, $\varphi \in \left( 0,\frac{\pi}{2} \right)$ is the angle
of deviation of the arms of the governor from its vertical axis
$S_1$, $\Omega \in [0,\infty)$ is the angular velocity of the
rotation of the engine flywheel $D$, $\theta$ is the angular
velocity of the rotation of  $S_1$, $l$ is the length of the arms,
$m$ is the mass of each ball, $H$ is a special sleeve, $T$ is a
set of transmission gears and $V$ is the valve that determines the
supply of steam to the engine.

The differential equations of our model which  generalize  those
found in Pontryagin  \cite{pon}, p. 217, are given by
\begin{eqnarray}\label{wattde}
\frac{d \; \varphi}{d \tau} &=& \psi \nonumber\\
\frac{d \; \psi}{d \tau} &=& (s(\Omega))^2 \; \sin \varphi \; \cos
\varphi -
\frac{g}{l} \; \sin \varphi - \frac{1}{m} \; h(\psi) \\
\frac{d \; \Omega}{d \tau} &=& \frac{1}{I} \; \left( M(\varphi) -F
\right) \nonumber
\end{eqnarray}
\noindent or equivalently by $ {\bf x}' = f ({\bf x})$,
\begin{equation}
f ({\bf x}) = \left( \psi, (s(\Omega))^2 \; \sin \varphi \; \cos
\varphi - \frac{g}{l} \; \sin \varphi - \frac{1}{m} \; h(\psi),
\frac{1}{I} \; \left( M(\varphi) - F \right) \right),
\label{campoinicial}
\end{equation}

where $\tau$ is the time, $\psi = \frac{d \varphi}{d \tau}$, $g$
is the standard acceleration of gravity, $\theta = s(\Omega) $,
$s$ is a smooth non-negative increasing function with $s(0) = 0 $,
called the {\it transmission function}, $h$ is a non-negative
increasing function with $h(0)=0$ that  represents the frictional
force of the system, $I$ is the moment of inertia of the flywheel
$D$ and  $M$ is a smooth decreasing function of the angle
$\varphi$. This function depends upon the design of the WGS and
determines the supply of steam to the engine through the valve
$V$. It measures the effect that the dynamics of the governor has
on the engine. Denny in \cite{denny} has also considered this
function for the torque. The torque $F$ due to the load is given
by $F = M (\varphi_\ast)$, where $\varphi_\ast$ is an equilibrium
angle.

The standard Watt governor differential equations in Pontryagin
\cite{pon}, p. 217,
\begin{eqnarray}\label{standard}
\frac{d \; \varphi}{d \tau} &=& \psi \nonumber\\
\frac{d \; \psi}{d \tau} &=& c^2 \; \Omega^2 \; \sin \varphi \;
\cos \varphi -
\frac{g}{l} \; \sin \varphi - \frac{b}{m} \; \psi \\
\frac{d \; \Omega}{d \tau} &=& \frac{1}{I} \; \left( \mu \cos
\varphi -F \right) \nonumber
\end{eqnarray}
are obtained from (\ref{wattde}) taking
\begin{equation}
s(\Omega)= c \; \Omega, \: h(\psi)= b \; \psi, \: M(\varphi) = \mu
\; \cos \varphi, \:\: c,b,\mu >0.
\end{equation}

\subsection{Stability analysis of the equilibrium
point}\label{stability}

Let $P_0 = (\varphi_0 , 0 , \Omega_0)$ be an equilibrium point of
(\ref{wattde}). The Jacobian matrix of $f$ at $P_0$ has the form
\begin{equation}
Df \left( P_0 \right) = \left( \begin{array}{ccc}
0 & 1 & 0 \\
\\- \displaystyle \frac{g \: \sin^2 \varphi_0}{l \: \cos \varphi_0} & -
\displaystyle \frac{h'(0)}{m} & \displaystyle \frac{2 \: g \:
s'(\Omega_0) \: \sin \varphi_0}{l \: s(\Omega_0)} \\
\\ \displaystyle \frac{M'(\varphi_0)}{I}    & 0  & 0
\end{array} \right).
\label{jacobianini}
\end{equation}

For the sake of completeness we state the following lemma whose
proof can be found in \cite{pon}, p. 58.

\begin{lema}
The polynomial $L(\lambda) = p_0 \lambda ^3 + p_1 \lambda ^2 + p_2
\lambda + p_3$, $p_0 > 0$, with real coefficients has all roots
with negative real parts if and only if the numbers $p_1 , p_2 ,
p_3$ are positive and the inequality $p_1 p_2 > p_0 p_3$ is
satisfied.

\label{routh}
\end{lema}

\begin{teo}
If
\begin{equation}
h'(0) > - \frac {2 \; m}{I} \: \frac{s'(\Omega_0)}{s(\Omega_0)} \:
M'(\varphi_0) \: \cot \varphi_0 \label{asymp}
\end{equation}
then the WGS differential equations (\ref{wattde}) have an
asymptotically stable equilibrium point at $P_0$. If
\[
0 < h'(0) < - \frac {2 \; m}{I} \:
\frac{s'(\Omega_0)}{s(\Omega_0)} \: M'(\varphi_0) \: \cot
\varphi_0
\]
then $P_0$ is unstable.

\label{teoestabilidade}
\end{teo}

\noindent {\bf Proof.} The characteristic polynomial of $Df \left(
P_0 \right)$ is given by $p(\lambda)$, where
\begin{equation}
-p(\lambda)= \lambda^3 + \frac{h'(0)}{m} \: \lambda^2 + \frac{g \:
\sin^2 \varphi_0}{l \: \cos \varphi_0} \: \lambda - \frac{2 \: g
\: M'(\varphi_0) \: s'(\Omega_0) \: \sin \varphi_0}{l \: I \:
s(\Omega_0)}. \label{charac}
\end{equation}
As $h'(0) > 0$, $s'(\Omega_0) > 0$, $ s (\Omega_0) > 0 $ and
$M'(\varphi_0) < 0$ the coefficients of $- p(\lambda)$ are
positive. Thus a necessary and sufficient condition for the
asymptotic stability of the equilibrium point $P_0$, as provided
by the condition for one real negative root and a pair of complex
conjugate roots with negative real part, is given by
(\ref{asymp}), according to Lemma \ref{routh}.
\begin{flushright}
$\blacksquare$
\end{flushright}
In terms of the WGS physical parameters, condition (\ref{asymp})
is equivalent to
\begin{equation}
\frac{h'(0) \: I}{m} \: \eta > 1, \label{asymparam}
\end{equation}
\noindent where
\[
\eta = \left| \frac{d \Omega_0}{dF} \right| = -
\frac{s(\Omega_0)}{2 \; M'(\varphi_0) \; s'(\Omega_0)} \: \tan
\varphi_0
\]
\noindent is the non-uniformity of the performance of the engine
which quantifies the change in the engine speed with respect to
the load (see \cite{pon}, p. 219, for more details). The rules
formulated by Vyshnegradskii to enhance the stability of the WGS
follow directly from (\ref{asymparam}). In particular, the
interpretation of (\ref{asymparam}) is that a sufficient amount of
damping $h'(0)$ must be present relative to the other physical
parameters for the system to be stable at the desired operating
speed. The general condition (\ref{asymparam}) is equivalent to
the original conditions given by Vyshnegradskii (see \cite{pon},
p. 219).

In next section we study the stability of $P_0$ under the
condition
\begin{equation}
h'(0) = - \frac {2 \; m}{I} \: \frac{s'(\Omega_0)}{s(\Omega_0)} \:
M'(\varphi_0) \: \cot \varphi_0, \label{valorcrit}
\end{equation}
\noindent that is, on the surface ---the Hopf surface---
complementary to the range of validity of Theorem
\ref{teoestabilidade}.

\section{Hopf bifurcation analysis}\label{hopf}

For the  analysis carried out here we take (\ref{wattde}) with
$s(\Omega) = c \; \Omega$, where  $c>0$ is a constant transmission
ratio and $h(\psi) = b \, \psi$, where  $b>0$ is a constant of the
frictional force of the system.

After the following change in the coordinates and the time
\begin{equation}
x = \varphi, \:\: y = \sqrt {\frac {l}{g}} \: \psi, \:\: z = c \;
\sqrt {\frac {l}{g}} \: \Omega, \:\: \tau = \sqrt {\frac {l}{g}}
\: t, \label{mudanca}
\end{equation}
\noindent the differential equations (\ref{wattde}) can be written
as
\begin{eqnarray}\label{wattdemm}
x' = \frac{d x}{d t} &=& y \nonumber\\
y' = \frac{d y}{d t} &=& z^2 \; \sin x \; \cos x -
\sin x - \varepsilon \; y \\
z' = \frac{d z}{d t} &=& T(x) - G \nonumber
\end{eqnarray}
\noindent or equivalently as ${\bf x}'= f({\bf x})$ where
\[
{\bf x} = (x,y,z) \in \left( 0, \frac{\pi}{2} \right) \times \R
\times [0,\infty),
\]
\begin{equation}
f({\bf x}) = \left(y, z^2 \; \sin x \; \cos x - \sin x -
\varepsilon \; y, T(x) - G \right), \label{campofinal}
\end{equation}
\[
\varepsilon = \frac {b}{m} \sqrt \frac{l}{g}, \: T(x) = \frac{c \;
l }{I g} \: M(x), \: G = \frac{c \; l \; F}{I g}.
\]
Here $\varepsilon$ is considered to be a changing  parameter; so
that the differential equations (\ref{wattdemm}) (or in its
equivalent vectorial form (\ref{campofinal}) ) can  in fact  be
regarded as a family of one-parameter families of differential
equations, dependent on the functional parameter $T$.

In this section we will analyze the stability at
\begin{equation}
P_0 = \left(x_0, y_0, z_0 \right)= \left( \arccos \beta, 0, \sqrt
\frac{1}{\beta}  \right),  \label{equilibriums}
\end{equation}
\noindent  under the condition (\ref{valorcrit}), which  now
writes as
\begin{equation}
\varepsilon_c = - \frac {2 \; \beta}{\omega_0} \: T'(x_0),
\label{epsiloncrit}
\end{equation}
where $ \beta = \cos x_0$ and
\begin{equation}
\omega_0 = \sqrt {\frac{1-\beta^2}{\beta}}. \label{omegazero}
\end{equation}

\subsection{Generalities on Hopf bifurcations}
The study  outlined below is based on the approach found in the
book of Kuznetsov \cite{kuznet}, pp 177-181.

Consider the differential equations
\begin{equation}
{\bf x}' = f ({\bf x}, {\bf \mu}), \label{diffequat}
\end{equation}
\noindent where ${\bf x} \in \R^3$ and ${\bf \mu} \in \R^m$ is a
vector of control parameters. Suppose (\ref{diffequat}) has an
equilibrium point ${\bf x} = {\bf x_0}$ at ${\bf \mu} = {\bf
\mu_0}$ and represent
\begin{equation}
F({\bf x}) = f ({\bf x}, {\bf \mu_0}) \label{Fhomo}
\end{equation}
as
\begin{equation}
F({\bf x}) = A{\bf x} + \frac{1}{2} \: B({\bf x},{\bf x}) +
\frac{1}{6} \: C({\bf x}, {\bf x}, {\bf x}) + O(|| {\bf x}
||^4){\nonumber}, \label{taylorexp}
\end{equation}
\noindent where $A = f_{\bf x}(0,{\bf \mu_0})$ and
\begin{equation}
B_i ({\bf x},{\bf y}) = \sum_{j,k=1}^3 \frac{\partial ^2
F_i(\xi)}{\partial \xi_j \: \partial \xi_k} \bigg|_{\xi=0} x_j \;
y_k, \label{Bap}
\end{equation}
\begin{equation}
C_i ({\bf x},{\bf y},{\bf z}) = \sum_{j,k,l=1}^3 \frac{\partial ^3
F_i(\xi)}{\partial \xi_j \: \partial \xi_k \: \partial \xi_l}
\bigg|_{\xi=0} x_j \; y_k \: z_l, \label{Cap}
\end{equation}
\noindent for $i = 1, 2, 3$. Here the variable ${\bf x}-{\bf x_0}$
is also denoted by ${\bf x}$.

Suppose $({\bf x_0}, {\bf \mu_0})$ is an equilibrium point of
(\ref{diffequat}) where the Jacobian matrix $A$ has a pair of
purely imaginary eigenvalues $\lambda_{2,3} = \pm i \omega_0$,
$\omega_0 > 0$, and no other critical (i.e., on the imaginary
axis)  eigenvalues.

The two dimensional center manifold can be parametrized by $w \in
\R^2 = \C$,  by means of ${\bf x} = H (w,\bar w )$, which is
written as
\[
H(w,{\bar w}) = w q + {\bar w}{\bar q} + \sum_{2 \leq j+k \leq 3}
\frac{1}{j!k!} \: h_{jk}w^j{\bar w}^k + O(|w|^4),
\]
with $h_{jk} \in \C ^3$, $h_{jk}={\bar h}_{kj}$.

Substituting these expressions into (\ref{diffequat}) and
(\ref{taylorexp}) we have
\begin{equation}
H_w (w,\bar w )w' + H_{\bar w}  (w,\bar w ){\bar w}'  = F(H(w
,\bar w )).
\label{homologicalp}
\end{equation}

Let $p, q \in \C ^3$ be vectors such that
\begin{equation}
A q = i \omega_0 \: q,\:\: A^{\top} p = -i \omega_0 \: p, \:\:
\langle p,q \rangle = \sum_{i=1}^3 \bar{p}_i \: q_i \:\: = 1.
\label{normalization}
\end{equation}

The complex vectors $h_{ij}$ are to be determined so that equation
(\ref{homologicalp}) writes as follows
\[
w'= i \omega_0 w + \frac{1}{2} \: G_{21} w |w|^2 + O(|w|^4),
\]
with  $G_{21} \in \C $.

Solving the  linear system  obtained by expanding
(\ref{homologicalp}),  the coefficients of the quadratic terms of
(\ref{Fhomo})  lead to
\begin{equation}
h_{11}=-A^{-1}B(q,{\bar q}) \label{h11},
\end{equation}
\begin{equation}
h_{20}=(2i\omega_0 I_3 - A)^{-1}B(q,q),\label{h20}
\end{equation}
where $I_3$ is the unit $3 \times 3$ matrix.

The coefficients  of the cubic terms are also uniquely calculated,
except  for  the term  $w^2 {\bar w}$, whose coefficient satisfies
a singular system for $h_{21}$
\begin{equation}
(i \omega_0 I_3 -A)h_{21}=C(q,q,{\bar q})+B({\bar q},h_{20}) + 2
B(q,h_{11})-G_{21}q, \label{h21m}
\end{equation}
which has a solution if and only if
\[
\langle p, C(q,q,\bar q) + B(\bar q, h_{20}) + 2 B(q,h_{11})
-G_{21} q \rangle = 0.
\]
Therefore
\begin{equation} \label{G21}
G_{21}= \langle p, C(q,q,\bar q) + B(\bar q, (2i \omega_0
I_3-A)^{-1} B(q,q)) - 2 B(q,A^{-1} B(q,\bar q)) \rangle,
\end{equation}
and the {\it first Lyapunov coefficient} $l_1$ -- which decides by
the analysis of third order terms at the equilibrium its
stability, if negative, or instability, if positive -- is defined
by
\begin{equation}
l_1 =  \frac{1}{2 \; \omega_0} \: {\rm Re} \; G_{21}.
\label{defcoef}
\end{equation}

A {\it Hopf point} $({\bf x_0}, {\bf \mu_0})$ is an equilibrium
point of (\ref{diffequat}) where the Jacobian matrix $A$ has a
pair of purely imaginary eigenvalues $\lambda_{2,3} = \pm i
\omega_0$, $\omega_0 > 0$, and no other critical eigenvalues. At a
Hopf point, a two dimensional center manifold  is well-defined,
which is invariant under the flow generated by (\ref{diffequat})
and can be smoothly continued to nearby parameter values.

A Hopf point is called {\it transversal} if the curves of complex
eigenvalues cross the imaginary axis with non-zero derivative.

In a neighborhood of a transversal Hopf point with $l_1 \neq 0$
the dynamic behavior of the system (\ref{diffequat}), reduced to
the family of parameter-dependent continuations of the center
manifold, is orbitally topologically equivalent to the complex
normal form
\begin{equation}\label{nf}
w' = (\gamma + i \omega) w + l_1 w |w|^2 ,
\end{equation}
$w \in \C $, $\gamma$, $\omega$ and $l_1$ are smooth continuations
of $0$, $\omega_0$ and the first Lyapunov coefficient at the Hopf
point \cite{kuznet}. When $l_1 < 0$ ($l_1
> 0$) a  family of stable (unstable) periodic orbits appears  found
on this family of center manifolds, shrinking to the  equilibrium
point at the Hopf point.

From (\ref{campofinal}) write the Taylor's expansion
(\ref{taylorexp}) of $f({\bf x})$. Define $a_1 = T'(x_0)$, $a_2 =
T''(x_0)$ and $a_3 = T'''(x_0)$. Thus
\begin{equation}
A = \left( \begin{array}{ccc}
0 & 1 & 0 \\
\\- \omega_0 ^2 & - \varepsilon_c & 2 \: \beta \: \omega_0 \\
\\ a_1    & 0  & 0
\end{array} \right),
\label{partelinear}
\end{equation}
\noindent and,  with the notation in (\ref{taylorexp}) we have
\begin{equation}
F({\bf x})\, - \, A{\bf x} = \left( 0,F_2({\bf x})+ O(||x||^4),
\frac{a_2}{2} \: x^2 + \frac{a_3}{6} \: x^3 + O(||x||^4) \right),
\label{partenaolinear}
\end{equation}
\noindent where
\begin{eqnarray*}
F_2({\bf x}) = - \frac{3}{2} \: \omega_0 \: \sqrt \beta \: x^2 +
\omega_0 \: \beta^{3/2} \: z^2 + \frac{2 (2 \beta ^2 -1)}{\sqrt
\beta} \: x \: z + \frac{4-7 \beta ^2}{6 \beta} \: x^3 - \\ 4 \:
\omega_0 \: \beta \: x^2 \: z  + (2 \beta ^2 -1) \: x \: z^2 .
\end{eqnarray*}

From (\ref{partelinear}) the eigenvalues of $A$ are
\begin{equation}
\lambda_1 = -\varepsilon_c , \: \: \lambda_2 = i \: \omega_0, \:
\: \lambda_3 = -i \: \omega_0. \label{autovalores}
\end{equation}
\noindent The eigenvectors $q$ and $p$ satisfying
(\ref{normalization}) are respectively
\begin{equation}
q = \left( -i, \omega_0, \frac{\varepsilon_c}{2 \beta} \right)
\label{q}
\end{equation}
\noindent and
\begin{equation}
p = \left( -\frac{i}{2}, \frac{\omega_0 - i \varepsilon_c}{2 (
\omega_0 ^2 + \varepsilon_c ^2 )}, \frac{\beta (\varepsilon_c + i
\omega_0)}{ \omega_0 ^2 + \varepsilon_c ^2 } \right). \label{p}
\end{equation}

The main result of this section can be formulated now.

\begin{teo}
Consider the family of differential equations (\ref{wattdemm}).
The first Lyapunov coefficient  at the point (\ref{equilibriums})
for parameter values satisfying (\ref{epsiloncrit}) is given by
\begin{equation}
l_1 (\beta, a_1, a_2, a_3, \varepsilon_c) = - \frac{R(\beta, a_1,
a_2, a_3)}{ 2 \left( {{(1- \beta^2)}^4} + 5 \beta^4 {{(1 -
\beta^2)}^2} a_1 ^2 + 4 \beta^8 a_1 ^4 \right)},
\label{coeficiente1}
\end{equation}
\noindent where
\begin{eqnarray*}
R(\beta, a_1, a_2, a_3) = \beta ^2 \Bigg[ 2 {{\beta }^7} a_1 ^6+
{{\beta }^{3}} a_1 ^4 \Big( 2 - 3 \beta^2 + 5
\beta^4 \Big) + 4 \beta {{(1-{{\beta }^2})}^3} a_2 ^2+ \\
\beta (1-{{\beta }^2}) \Big( a_1 ^2 \Big( 9 (1-{{\beta }^2})^2 + 2
\beta^4 a_2 ^2 \Big) + 2  a_3 \Big( - \beta^4 a_1 ^3 - (1 -
\beta^2)^2 a_1 \Big) \Big) + \\
a_2 \; \sqrt {1 - \beta^2 } \Big( - \beta^4 a_1 ^3 (1+ 5 \beta^2)-
a_1 (1 - \beta^2)^2 (5 + 3 \beta^2) \Big)  \Bigg].
\end{eqnarray*}

\label{lemacoef1}
\end{teo}

\noindent{\bf Proof.} The proof depends on preliminary
calculations presented below. From (\ref{taylorexp}), (\ref{Bap}),
(\ref{Cap}) and (\ref{partenaolinear}) one has
\begin{equation}
B({\bf x},{\bf y}) = \left( 0, B_2({\bf x},{\bf y}), a_2 \: x_1 \:
y_1 \right), \label{B1}
\end{equation}
\noindent where
\begin{eqnarray*}
B_2({\bf x},{\bf y})= -3 \: \omega_0 \: \sqrt \beta \: x_1 \: y_1
+ 2 \: \omega_0 \: \beta^{3/2} \: x_3 \: y_3 + \frac{2 \: (2 \:
\beta ^2 -1)}{\sqrt \beta} \left( x_1 \: y_3 + x_3 \: y_1 \right),
\end{eqnarray*}
\begin{equation}
C({\bf x},{\bf y}, {\bf z}) = \left( 0, C_2({\bf x},{\bf y},{\bf
z}), a_3 \: x_1 \: y_1 \: z_1 \right), \label{C1}
\end{equation}
where
\begin{eqnarray*}
C_2({\bf x},{\bf y},{\bf z})= \frac{4-7 \beta^2}{\beta} \: x_1 \:
y_1 \: z_1 - 8 \omega_0 \: \beta \left( x_1 \: y_1 \: z_3 + x_1 \:
y_3 \: z_1 + x_3 \: y_1 \: z_1 \right) + \\ 2 \: (2 \beta^2 -1)
\left(x_1 \: y_3 \: z_3 + x_3 \: y_1 \: z_3 + x_3 \: y_3 \: z_1
\right).
\end{eqnarray*}

Referring to the notation in (\ref{B1}), (\ref{C1}), (\ref{h11}),
(\ref{h20}) and (\ref{q}) one has
\begin{equation}
B(q,q) = \left( 0, \displaystyle \frac{\beta \omega_0
(\varepsilon_c ^2 + 6 \beta) - i 4 \varepsilon_c (2 \beta^2 -1)}{2
\beta^{3/2}}, -a_2 \right), \label{Bqq}
\end{equation}
\begin{equation}
B(q,\bar{q}) = \left( 0, \displaystyle \frac{\omega_0
(\varepsilon_c ^2 -6 \beta)} {2 \sqrt \beta},  a_2 \right),
\label{Bqqbar}
\end{equation}
\begin{equation}
C(q,q,\bar{q}) = \left( 0,  \displaystyle \frac{ - 8 \beta^2
\varepsilon_c \omega_0 + i \left( \beta (14 \beta^2 - 8) -
\varepsilon_c ^2 (2 \beta^2 -1) \right)}{2 \beta^2}, -i \; a_3
\right), \label{Cqqqbar}
\end{equation}
\begin{equation}
B(q,h_{11}) = \left( 0, B_2(q,h_{11}), i \: \frac{a_2 ^2}{a_1}
\right), \label{Bqr}
\end{equation}
\noindent where
\begin{eqnarray*}
B_2 (q,h_{11})= \frac{1}{4 \beta^2 a_1} \Bigg[(6 \beta -
\varepsilon_c ^2) (2 i -4 i \beta^2 + \beta \varepsilon_c
\omega_0) a_1 - 2 {\sqrt \beta} \Bigg( 2 i (1+\beta^2) \omega_0 +\\
\varepsilon_c \left( -2 + \beta (4 \beta + \omega_0 ^2) \right)
\Bigg) a_2 \Bigg],
\end{eqnarray*}
\noindent and
\begin{equation}
B(\bar{q},h_{20}) = \left(  0, B_2(\bar{q},h_{20}),
B_3(\bar{q},h_{20}) \right), \label{Bqbars}
\end{equation}
\noindent where
\begin{eqnarray*}
B_2(\bar{q},h_{20})= \frac{1}{4 \beta^3 \omega_0 ( (2
\varepsilon_c + 3 i \omega_0) \omega_0 + \beta a_1)} \Bigg[ \beta
\Bigg( 12 i \beta^2 (2 \beta^2 -1) \omega_0 + \\ 6 i \beta (2
\beta^2 -1) \varepsilon_c ^2 \omega_0 - \beta^2 \varepsilon_c ^3
\omega_0 ^2 + \varepsilon_c \left( 8 - 32 \beta^2 + 32 \beta^4 - 6
\beta^3 \omega_0 ^2 \right) \Bigg)a_1 + \\ 2 \omega_0 \Bigg( -18 i
\beta^2 (2 \beta^2 -1) \varepsilon_c \omega_0 - i \beta (2 \beta^2
-1) \varepsilon_c ^3 \omega_0 +  18 \beta^4 \omega_0 ^2 + \\
\varepsilon_c ^2 \left( -4 + 16 \beta^2 -16 \beta^4 + 3 \beta^3
\omega_0 ^2 \right)+ \beta^{5/2} \Bigg( 6 i
(3 \beta^2 -1) \omega_0 + 2 i \beta \omega_0 \varepsilon_c ^2 + \\
3 \varepsilon_c (-2 + 4 \beta^2 - \beta \omega_0 ^2) \Bigg) a_2
\Bigg) \Bigg],
\end{eqnarray*}
\noindent and
\[
B_3(\bar{q},h_{20}) = - \frac{ a_2 \left( \beta \omega_0 (6 \beta
+ \varepsilon_c ^2) + i \left( (4-8 \beta^2) \varepsilon_c + 2
\beta ^{5/2} a_2 \right) \right) }{2 \beta^{3/2} \left(2
\varepsilon_c \omega_0 + \beta a_1 + 3 i \omega_0 ^2 \right)}.
\]

The first Lyapunov coefficient is given by (\ref{defcoef}). From
(\ref{p}) and (\ref{Cqqqbar}) one has
\begin{equation}
{\rm Re} \langle p, C(q,q,\bar{q}) \rangle = \frac{\varepsilon_c
\left( 8 \beta - 14 \beta^3 - \varepsilon_c ^2 + 2 \beta^2
(\varepsilon_c ^2 - 4 \omega_0 ^2) \right) - 4 \beta^3 \omega_0
a_3}{4 \beta^2 (\omega_0 ^2 + \varepsilon_c ^2)}.
\label{partereal1}
\end{equation}
\noindent From (\ref{p}) and (\ref{Bqr}) one has
\begin{eqnarray}\label{partereal2}
{\rm Re} \langle p, 2 B(q,h_{11}) \rangle = \frac{1}{4 \beta^2
(\omega_0 ^2 + \varepsilon_c ^2) a_1} \Bigg[ \varepsilon_c (6
\beta -\varepsilon_c ^2) (-2 + 4 \beta^2 + \beta \omega_0 ^2) a_1
+ \nonumber \\ 2 \sqrt{\beta} \omega_0 a_2 \left( \varepsilon_c (4
-2 \beta^2 -\beta \omega_0 ^2) + 4 \beta^{5/2} a_2 \right) \Bigg].
\end{eqnarray}
\noindent From   (\ref{p}) and (\ref{Bqbars}) one has
\begin{equation}
{\rm Re} \langle p, B(\bar{q},h_{20}) \rangle =
\frac{\vartheta(\beta, a_1, a_2, \varepsilon_c, \omega_0)}{8
\beta^3 \left( \omega_0 ^2 + \varepsilon_c ^2 \right) \left( 4
\varepsilon_c \omega_0 ^2 + 9 \omega_0 ^4 + \beta a_1 (4
\varepsilon_c \omega_0 + \beta a_1) \right)}, \label{partereal3}
\end{equation}
where
\begin{eqnarray*}
\vartheta(\beta, a_1, a_2 , \varepsilon_c, \omega_0)= 4 (2 \beta^2
-1) \varepsilon_c ^3 \left( 10 + \beta (\varepsilon_c ^2 - 26
\beta) \right) \omega_0 ^2 + \beta \Bigg( a_1 \Bigg( 2 (1 - \beta
^2) \\ \varepsilon_c ^2 ( \beta (42 \beta + \varepsilon_c ^2 ) -
24) \omega_0 - \beta (36 \beta (3 \beta ^2 -1) + 6 (1 + 4 \beta
^2) \varepsilon_c ^2 + 5 \beta \varepsilon_c ^4)\omega_0 ^3 - \\
\beta \varepsilon_c (2 (2 \beta ^2 -1) (\beta (2 \beta +
\varepsilon_c ^2) - 4) +  \beta ^2 (6 \beta + \varepsilon_c ^2)
\omega_0 ^2)a_1 \Bigg) + 2 \beta ^{3/2} \omega_0 \\ \left( 18 (1 -
5 \beta ^2) \omega_0 ^3 + \varepsilon_c ^2 \omega_0 (38 - 5 \beta
(8 \beta + 3 \omega_0 ^3)) + \beta \varepsilon_c (4 + 10 \beta ^2
- 3 \beta \omega_0 ^2) a_1 \right) a_2 + \\ 8 \beta ^4 \omega_0 (5
\varepsilon_c \omega_0 + \beta a_1) a_2 ^2 \Bigg) - 6 \beta
\varepsilon_c \omega_0 ^4 \left( 6 \beta (1 + 3 \beta ^2) + (7
\beta ^2 -1) \varepsilon_c ^2 \right).
\end{eqnarray*}
\noindent Substituting (\ref{omegazero}) and (\ref{epsiloncrit})
into (\ref{partereal1}), (\ref{partereal2}) and (\ref{partereal3})
and the results into (\ref{G21}) and (\ref{defcoef}), the theorem
is proved.
\begin{flushright}
$\blacksquare$
\end{flushright}

\begin{prop}
Consider the family of differential equations (\ref{wattdemm})
regarded as dependent on the parameter $\varepsilon$. The real
part, $\gamma$, of the pair of complex  eigenvalues  verifies
\begin{equation}
\gamma'(\varepsilon_c) = -\frac{1}{2} \frac{\omega_0 ^2}{\omega_0
^2 + \varepsilon_c ^2} < 0. \label{transversal}
\end{equation}
Therefore, the transversality condition holds at the Hopf point.
\label{lematransv}
\end{prop}

\begin{figure}[!h]

\centerline{\includegraphics[width=8cm]{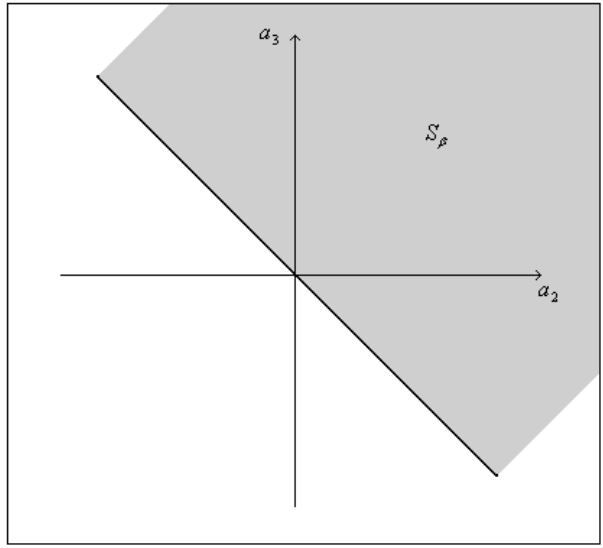}}

\caption{{\small The region $S_{\beta}$ }.}

\label{Sbeta}

\end{figure}

\noindent{\bf Proof.} Let $\lambda (\varepsilon) = \lambda_{2,3}
(\varepsilon) = \gamma (\varepsilon) \pm i \omega (\varepsilon)$
be eigenvalues of $A(\varepsilon)$ such that $\gamma
(\varepsilon_c) = 0$ and $\omega (\varepsilon_c) = \omega_0$,
according to (\ref{autovalores}). Taking the inner product of $p$
with the derivative of
\[
A(\varepsilon) q(\varepsilon) = \lambda (\varepsilon)
q(\varepsilon)
\]
at $\varepsilon = \varepsilon_c$ one has
\[
\left \langle p , \frac{d A}{d \varepsilon} \Bigg |_ {\varepsilon
= \varepsilon_c} \: q \right \rangle = \gamma'(\varepsilon_c) \pm
\omega'(\varepsilon_c).
\]
Thus the transversality condition is given by
\begin{equation}
\gamma'(\varepsilon_c) = {\rm Re} \: \left \langle p, \frac{d A}{d
\varepsilon} \Bigg |_ {\varepsilon = \varepsilon_c} \: q \right
\rangle. \label{deftransv}
\end{equation}
\noindent As
\[
\frac {d A}{d \varepsilon} \Bigg |_ {\varepsilon = \varepsilon_c}
\: q = \left( 0, - \omega_0, 0 \right),
\]
the proposition follows from a simple calculation.
\begin{flushright}
$\blacksquare$
\end{flushright}
\begin{teo}
Define $S_{\beta} = \{(a_2,a_3)| \: a_3 + K_2 (\beta) a_2
> 0 \}$, where
\[
K_2(\beta) = \frac{5+3 \beta^2}{2 \beta \sqrt{1 - \beta^2}}.
\]
If $(a_2,a_3) \in S_{\beta}$ then the one-parameter family of
differential equations (\ref{wattdemm}) has a Hopf point at $P_0$
for $\varepsilon = \varepsilon_c$. Furthermore this Hopf point at
$P_0$ is asymptotically stable and for each $\varepsilon <
\varepsilon_c$, but close to $\varepsilon_c$, there exists a
stable periodic orbit near the unstable equilibrium point $P_0$.
See Fig. \ref{Sbeta}.

\label{teohopf}
\end{teo}

\noindent{\bf Proof.} Since the denominator of
(\ref{coeficiente1}) is positive the sign of the first Lyapunov
coefficient is determined by the sign of $R(\beta, a_1, a_2,
a_3)$. Rewrite this expression in the following way
\begin{eqnarray*}
\frac{R(\beta, a_1, a_2, a_3)}{\beta^2} = \Bigg[ 2 {{\beta }^7}
a_1 ^6+ {{\beta }^{3}} a_1 ^4 \Big( 2 - 3 \beta^2 + 5
\beta^4 \Big) + 4 \beta {{(1-{{\beta }^2})}^3} a_2 ^2+ \\
\beta (1-{{\beta }^2}) \Big( a_1 ^2 \Big( 9 (1-{{\beta }^2})^2 + 2
\beta^4 a_2 ^2 \Big) \Big)  \Bigg] - 2 \beta^5 a_1 ^3 (1 -
\beta^2) (a_3 + K_1(\beta) a_2)- \\ 2 a_1 \beta (1-\beta^2)^3 (a_3
+ K_2(\beta) a_2),
\end{eqnarray*}
where
\[
K_1 (\beta) = \frac{1 + 5 \beta^2}{2 \beta \sqrt{1 - \beta^2}} >0,
\]
and $K_2(\beta)$ is as above. Since $ 2 - 3 \beta^2 + 5 \beta^4$
is positive, the sum into the bracket is positive. A simple
calculation shows that
\[
K_2(\beta)> K_1(\beta),
\]
for all $\beta \in (0,1)$. If $(a_2,a_3) \in S_{\beta}$ then
\[
- 2 \beta^5 a_1 ^3 (1 - \beta^2) (a_3 + K_1(\beta) a_2)- \\ 2 a_1
\beta (1-\beta^2)^3 (a_3 + K_2(\beta) a_2) > 0,
\]
since $a_1 < 0$. This implies that $R(\beta, a_1, a_2, a_3) > 0 $
and therefore the first Lyapunov coefficient is negative. See Fig.
\ref{Sbeta}.
\begin{flushright}
$\blacksquare$
\end{flushright}

\begin{remark}\label{pondenny}
The family (\ref{wattdemm}) has been studied in Denny
\cite{denny}, focusing the Lyapunov stability of the equilibrium
point (\ref{equilibriums}), outside the bifurcation surface
(\ref{epsiloncrit}). The stability conditions in \cite{pon} and
\cite{denny} have been extended in Theorem \ref{teoestabilidade}
to a more general WGS model given in (\ref{wattde}). For this more
general model, however, the study of the codimension one Hopf
bifurcation involves much longer calculations. For this reason,
the analysis performed here has been restricted  to the case
proposed by \cite{denny}. See Theorem \ref{teohopf}. The
calculations for the general case in (\ref{wattde}), while too
long to put in print, can be handled by Computer Algebra. The
special case of the equations of Pontryagin \cite{pon}  will be
treated in next subsection.
\end{remark}

\subsection{Hopf Bifurcations in Pontryagin equations}

Below we specialize  the calculations above  to the case of the
Eq. (\ref{standard}) for the WGS  in the form presented by
Pontryagin \cite{pon} also treated in \cite{humadi}.

With the change in the coordinates and time (\ref{mudanca}) the
standard Watt governor differential equations (\ref{standard}) are
obtained from (\ref{wattdemm}) with $M(x) = \mu \cos x$ and have
the form
\begin{eqnarray}\label{standardm}
x' &=& y \nonumber\\
y' &=& z^2 \; \sin x \; \cos x -
\sin x - \varepsilon \; y \\
z' &=& \alpha \; \left( \cos x - \beta \right) \nonumber
\end{eqnarray}
\noindent where $\alpha > 0$, $\beta \in (0,1)$ and $\varepsilon
> 0 $ are given by
\begin{equation}
\varepsilon = \frac{b}{m} \; \sqrt \frac{l}{g}, \:\: \alpha =
\frac{c \; l \; \mu}{g \; I}, \:\: \beta = \frac{F}{\mu}.
\label{newparameters}
\end{equation}

The differential equations (\ref{standardm}) have an equilibrium
point  located at $P_0$ given by Eq. (\ref{equilibriums}).

A necessary and sufficient condition for the hyperbolic asymptotic
stability of the equilibrium point $P_0$ is $\varepsilon > 2 \:
\alpha \: \beta ^{3/2}$, according to Theorem
\ref{teoestabilidade}. For $0 < \varepsilon < 2 \: \alpha \: \beta
^{3/2}$, the equilibrium point $P_0$ is unstable. Now we analyze
the stability of $P_0$ as $\varepsilon_c = \varepsilon
(\beta,\alpha) = 2 \: \alpha \: \beta ^{3/2}$ obtained from
(\ref{epsiloncrit}).

\begin{teo}
Consider the family of differential equations (\ref{standardm}).
The first Lyapunov coefficient on the Hopf surface $\varepsilon
=\varepsilon_c = 2\: \alpha \: \beta^{3/2}$ is given by
\begin{equation}
l_1 (\beta,\alpha,\varepsilon_c) = -  \frac{\alpha \beta^2 \sqrt{1
- \beta ^2} \left( 3 + (\alpha^2 -5) \beta ^2 + \alpha ^4 \beta ^6
\right)}{2 \left( 1 - \beta^2 + \alpha^2 \beta^4 \right) \left( 1
- \beta^2 + 4 \alpha^2 \beta^4 \right)}. \label{coefstandard}
\end{equation}

\label{l1standard}
\end{teo}

\noindent{\bf Proof.} The proof follows from (\ref{coeficiente1})
with
\[
a_1 = - \alpha \sqrt{1 - \beta^2}, \: \: a_2 = - \alpha \; \beta,
\: \: a_3 = \alpha \sqrt{1 - \beta^2}.
\]
\begin{flushright}
$\blacksquare$
\end{flushright}
\begin{teo}
If
\begin{equation}
g(\beta, \alpha) =  3 + (\alpha^2 -5) \beta ^2 + \alpha ^4 \beta
^6 \label{hopfcondstan}
\end{equation}
\noindent is different from zero then the family of differential
equations (\ref{standardm}) has a Hopf point at $P_0$ for
$\varepsilon_c = 2 \: \alpha \: \beta ^{3/2}$.

\label{curvagstan}
\end{teo}

\noindent {\bf Proof.} From (\ref{transversal}) the transversality
condition is satisfied. Therefore a sufficient condition for being
a Hopf point is that the first Lyapunov coefficient
$l_1(\beta,\alpha, \varepsilon_c) \neq 0$. But from
(\ref{coefstandard}) it is equivalent to $g(\beta, \alpha) \neq
0$. The theorem is proved.
\begin{flushright}
$\blacksquare$
\end{flushright}

The following theorem summarizes the results of this subsection.

\begin{teo}
If $(\beta, \alpha, \varepsilon_c) \in S \cup U$ then the family
of differential equations (\ref{standardm}) has a Hopf point at
$P_0$. If $(\beta, \alpha, \varepsilon_c) \in S$ then the Hopf
point at $P_0$ is asymptotically stable and for each $\varepsilon
< \varepsilon_c$, but close to $\varepsilon_c$, there exists a
stable periodic orbit near the unstable equilibrium point $P_0$.
If $(\beta, \alpha, \varepsilon_c) \in U$ then the Hopf point at
$P_0$ is unstable and for each $\varepsilon > \varepsilon_c$, but
close to $\varepsilon_c$, there exists an unstable periodic orbit
near the asymptotically stable equilibrium point $P_0$.
\label{teohopfstan}
\end{teo}

\begin{cor}
Consider the family of differential equations (\ref{standardm}).
If $\alpha > 1$ then the equilibrium point $P_0$ is asymptotically
stable for $\varepsilon = \varepsilon_c$ and for all $0 < \beta <
1$. Therefore for each $\varepsilon < \varepsilon_c$, but close to
$\varepsilon_c$, there exists a stable periodic orbit near the
unstable equilibrium point $P_0$.

\label{alfamaiorque1stan}
\end{cor}

\noindent {\bf Proof.} The proof is immediate from Theorem
\ref{teohopfstan}. In fact, if $\alpha > 1$ then $(\beta, \alpha,
\varepsilon_c) \in S$ for all $0 < \beta <1$. Thus the first
Lyapunov coefficient $l_1(\beta, \alpha, \varepsilon_c)$ is
negative.
\begin{flushright}
$\blacksquare$
\end{flushright}
\begin{figure}[!h]
\centerline{\includegraphics[width=9cm]{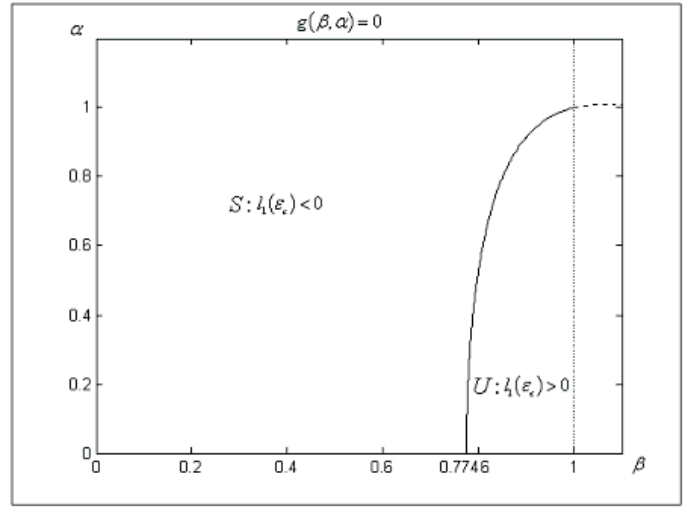}}
\caption{{\small Signs of the first Lyapunov coefficient}.}

\label{curvag=0}

\end{figure}
\begin{remark}\label{crucial}
Equation (\ref{hopfcondstan}) gives a simple expression to
determine  the sign of the first Lyapunov coefficient
(\ref{coefstandard}). Its graph is illustrated in Fig.
\ref{curvag=0}, where the signs of the first Lyapunov coefficient
are also represented. This gives an analytic  corroboration to the
observation in \cite{has1}, p. 155, made  on the basis of
numerical testing:

``For $ \alpha
>1 $, our computations indicate that $ l_1(\beta, \alpha,
\varepsilon_c) < 0 $ for all values of $ \beta $, $ 0 < \beta < 1
$."

This matter has been addressed in \cite{humadi}, providing a
longer expression to decide the sign first Lyapunov coefficient
than the one given in (\ref{coefstandard}).

The curve $l_1 = 0$ divides the surface of critical parameters
$\varepsilon_c = 2 \: \alpha \: \beta ^{3/2}$ into two connected
components denoted by $S$ and $U$ where $l_1 < 0$ and $l_1 > 0$
respectively.  See Fig. \ref{curvag=0}.
\end{remark}

\subsection{Examples of stability  and instability in the WGS} \label{specific}

The following examples give a piece of information about the
stability of the bifurcating periodic orbit of the family
(\ref{standardm}). First in terms of the angle of deviation
$\varphi$.
\begin{example}
If the equilibrium angle of deviation
\begin{equation}
\varphi_0 > 39.23^\circ \label{angulocriticostan}
\end{equation}
\noindent then the equilibrium point $P_0$ is asymptotically
stable for $\varepsilon = \varepsilon_c$ and for all $\alpha > 0$.
Therefore for each $\varepsilon < \varepsilon_c$, but close to
$\varepsilon_c$, there exists a stable periodic orbit near the
unstable equilibrium point $P_0$.
\label{alfaqualquerstan}
\end{example}
\noindent In fact, if $\varphi_0 > 39.23^\circ$ then $\beta = \cos
\varphi_0 < 0.7746$. Therefore $l_1 (\varepsilon_c) < 0$ for all
$\alpha > 0$,  as follows from  Fig. \ref{curvag=0}.

Now in terms of  the  normalized speed.
\begin{example}
If the normalized engine speed $z_0$ satisfies
\begin{equation}
z_0 > 1.1362 \label{enginespeed}
\end{equation}
\noindent then the equilibrium point $P_0$ is asymptotically
stable for $\varepsilon = \varepsilon_c$. Therefore for each
$\varepsilon < \varepsilon_c$, but close to $\varepsilon_c$, there
exists a stable periodic orbit near the unstable equilibrium point
$P_0$. \label{velocidadenormalstan}
\end{example}
\noindent In fact, if  $z_0 > 1.1362$ then $\varphi_0 >
39.23^\circ$ since $\cos \varphi_0 = 1/{z_0 ^2}$. The proof
follows from the Example \ref{alfaqualquerstan}.

In terms of the normalized parameter $\beta$  in
(\ref{newparameters}), we have:
\begin{example}
If
\[
0.7746 < \beta < 1
\]
\noindent and
\[
0 < \alpha < \frac{\sqrt { \sqrt {20 \beta^4 -12 \beta^2 + 1}
-1}}{\sqrt 2 \beta^2}
\]
\noindent then the equilibrium point $P_0$ is unstable for
$\varepsilon = \varepsilon_c$. Therefore for each $\varepsilon >
\varepsilon_c$, but close to $\varepsilon_c$, there exists an
unstable periodic orbit near the asymptotically stable equilibrium
point $P_0$. \label{equinstavelstan}
\end{example}
\noindent This follows from Fig. \ref{curvag=0}. In fact, the
level curve $g(\beta, \alpha) = 0$, where $g (\beta, \alpha) $ is
defined in (\ref{hopfcondstan}), is the graph of the function
\[ h(\beta) = \frac{\sqrt { \sqrt {20 \beta^4 -12 \beta^2
+ 1} -1}}{\sqrt 2 \beta^2}.
\]

\section{Concluding  comments}\label{conclusion}

The historical relevance of the Watt governor study as well as its
importance for present day theoretical and technological control
developments -- going from steam to diesel and gasoline engines --
have been widely discussed by MacFarlane \cite{mac}, Denny
\cite{denny} and Wellstead - Readman \cite{control} among others.

In this paper the original stability analysis due to Maxwell and
Vyshnegradskii of the Watt Centrifugal Governor System --WGS-- has
been revisited and extended to the following situations:

\begin{enumerate}
\item A more general torque function of the engine has been
considered. This function determines the supply of steam to the
engine via the valve $V$ and is related to the design of the WGS.
It measures the effect of the governor on the engine;
 \item A more general
transmission function has been  studied. This function determines
the relation between the angular velocities of the engine and the
governor axis. It measures the effect of the engine on the
governor;
 \item A more general frictional force of the system has been
modelled.
\end{enumerate}

In Theorem \ref{teoestabilidade} we have  extended  the stability
results presented  in Pontryagin \cite{pon} and Denny
\cite{denny}.

Concerning the bifurcations of the WGS  this paper deals with the
codimension one Hopf bifurcations in the Watt governor
differential equations. The main results are:
\begin{enumerate}
\item In Theorem \ref{teohopf} we give sufficient conditions for
the stability of the periodic orbit that appears from the Hopf
point for the general Watt governor differential equations
(\ref{wattdemm});

\item In Theorem \ref{teohopfstan} we have extended and provided a
neat geometric synthesis of the Hopf stability analysis performed
by Hassard et al. \cite{has1} and Al-Humadi and Kazarinoff
\cite{humadi}. This identifies the codimension one Hopf points and
the location of its complement codimension 2 Hopf points along a
curve in Fig. \ref{curvag=0}.
\end{enumerate}

\vspace{0.2cm} \noindent {\bf Acknowledgement}: The first and
second authors developed this work under the project CNPq Grant
473824/04-3. The first author is fellow of CNPq and takes part in
the project CNPq PADCT 620029/2004-8. This work was finished while
he visited Brown University, supported by FAPESP, Grant
05/56740-6.

\vspace{0.5cm}

\noindent Jorge Sotomayor

\noindent{\em Instituto de
Matem\'atica e Estat\'{\i}stica, Universidade de S\~ao Paulo\\ Rua
do Mat\~ao 1010, Cidade Universit\'aria\\ CEP 05.508-090, S\~ao
Paulo, SP, Brazil
\\}e--mail:sotp@ime.usp.br

\vspace{0.5cm}

\noindent Luis Fernando Mello

\noindent{\em Instituto de Ci\^encias Exatas \\Universidade
Federal de Itajub\'a \\CEP 37.500-903, Itajub\'a, MG, Brazil
\\}e--mail:lfmelo@unifei.edu.br

\vspace{0.5cm}

\noindent Denis de Carvalho Braga

\noindent{\em Instituto de Sistemas El\'etricos e Energia \\
Universidade Federal de Itajub\'a
\\CEP 37.500-903, Itajub\'a, MG, Brazil
\\}e--mail:braga\_denis@yahoo.com.br

\end{document}